\documentclass[11pt]{article}
\usepackage[a4paper]{geometry}
\usepackage{times,amsmath,amssymb,amsfonts,amsthm}

\makeatletter
\numberwithin{equation}{section} 
\numberwithin{figure}{section} 
\theoremstyle{plain}
\newtheorem{thm}{Theorem}
\theoremstyle{remark}
\newtheorem{rem}[thm]{Remark}
\newcommand{\IntO}{\int_{\Omega}}
\newcommand{\be}{\begin{eqnarray*}}
\newcommand{\ee}{\end{eqnarray*}}
\newcommand{\ben}{\begin{eqnarray}}
\newcommand{\een}{\end{eqnarray}}

\newcommand{\mV}{\mathcal{V}} 

\newcommand{\Div}{\mathrm{div}} 
\newcommand{\Hdiv}[1]{H_{\Div}({#1})} 
\newcommand{\R}{\mathbb{R}} 

\newcommand{\LN}[1]{\left\Vert{#1}\right\Vert} 
\newcommand{\CN}[1]{\left\vert\!\left\vert\!\left\vert{#1}\right\vert\!\right\vert\!\right\vert} 

\newcommand{\mM}{\mathcal{M}} 
\newcommand{\mN}{\mathcal{N}} 

\makeatother

\begin{document}

\title{A Posteriori Error Estimates for Nonconforming Approximations of Evolutionary Convection-Diffusion Problems}

\author{\textbf{S.I.~Repin} \\
V.A.~Steklov Institute of Mathematics,\\
Fontanka 27, 191023 St. Petersburg, Russia\\
repin@pdmi.ras.ru\\[2ex]
\textbf{S.K.~Tomar \footnote{Corresponding author}} \\
RICAM, Austrian Academy of Sciences,\\
Altenbergerstr. 69, 4040 Linz, Austria\\
satyendra.tomar@ricam.oeaw.ac.at}

\date{}

\maketitle

\begin{abstract}
We derive computable upper bounds for the difference between an exact solution of the evolutionary convection-diffusion problem and an approximation of this solution. The estimates are obtained by certain transformations of the integral identity that defines the generalized solution. These estimates depend on neither special properties of the exact solution nor its approximation, and involve only global constants coming from embedding inequalities.
The estimates are first derived for functions in the corresponding energy space, and then possible extensions to classes of piecewise continuous approximations are discussed.

\noindent
\textbf{Keywords:} a posteriori error estimates, time-nonconforming approximations, convection-diffusion problems

\end{abstract}

\section{Introduction}
Let $u(x,t)$ be the generalized solution of the initial-boundary value problem
\begin{equation}
u_t={\mathcal L} u \quad \mathrm{in } \,\, Q_T, \qquad u(x,0)=\varphi(x),
\end{equation}
where $Q_{T}:=\left(0,T\right) \times \Omega$, $\Omega$ is an open bounded domain in $\R^{d}$ with the Lipschitz boundary $\partial \Omega$, ${\mathcal L}$ is a linear uniformly elliptic operator, and $\varphi$ is a certain function that defines the initial condition. For this class of evolutionary problems the existence and regularity theory has been deeply elaborated (see, e.g., \cite{Ladyzhen-85, LadyzhenSU-67}). It is well known that if the boundary conditions on $S_{T}:=\left(0,T\right) \times \partial \Omega$ are defined in a suitable way (e.g., $u$ vanishes on $S_{T}$), then the generalized solution exists and is unique. Moreover, under some conditions imposed on the coefficients of ${\mathcal L}$, $\varphi$, and $\partial\Omega$, one can establish the so-called first and second main inequalities that estimate different norms of $u$.
In this paper, we focus on a different problem related to evolutionary models in mathematical physics. Assume that $v$ is a given function (e.g., a numerical solution, or a solution of some simplified mathematical model associated with the same phenomenon) which we wish to compare with $u$. In general, the conditions that we impose on $v$ are rather broad: we assume that it can be any function from the corresponding energy space. We do not suppose that $v$ is subject to some extra regularity assumptions or that it must satisfy Galerkin orthogonality property with respect to some finite-dimensional space. Estimates of such a type are derived by purely functional methods without using specific properties of approximations and do not involve mesh dependent constants. For this reason, they are sometimes called \textit{functional type a posteriori estimates}.
Finding computable and efficient estimates of a certain measure of $u-v$ presents a mathematical problem of high importance for quantitative analysis of evolutionary problems in mathematical physics. We note that for linear and nonlinear elliptic problems (including stationary reaction-convection-diffusion problems and variational inequalities) this problem is well studied (cf. a systematic exposition in \cite{Repin-08} and papers cited therein). In \cite{Repin-02}, it was shown that transformations of the integral identity that defines the generalized solution of the parabolic heat equation lead to directly computable bounds for some weighted space-time norms of $u-v$. In \cite{GaevsRepin-05}, such estimates were obtained (and numerically tested) for a wider class of linear parabolic problems. Recently, in \cite{NeittanRepin-10}, similar estimates were derived for the evolutionary Stokes problem.
In this paper, we consider the evolutionary convection-diffusion problem. We apply the same method as in \cite{Repin-02}, and show that for any $v$ in the energy function space computable and guaranteed upper bounds of certain norms of $u-v$ follow from the corresponding integral identity that defines the generalized solution. Moreover, we show that $v$ may not exactly satisfy the initial boundary conditions (this matter is discussed in the last section of the paper) and it may be discontinuous in time.

We consider the initial boundary value problem
\begin{alignat}{2}
u_{t} - \Div\nabla u + a \cdot \nabla u & =f & & \quad\mbox{in } Q_{T},\label{eq:Prob}\\
u\left(x,0\right) & =\varphi\left(x\right) &  & \quad x\in\Omega,\label{eq:Cond_Init}\\
u & =0 & & \quad\mbox{on } S_{T}, \label{eq:Cond_Bdry}
\end{alignat}
where $\Omega$ is an open bounded domain in $\R^{d}$. We also assume that
\begin{eqnarray}
&&a=a\left(x\right)\in L^{\infty}\left(\Omega,\R^{d}\right),\,
\Div a\in L^{\infty}\left(\Omega\right),\,
\Div a\le 0,\\
&& f \in L^{2,1}(Q_T), \qquad \varphi \in H^{1}_{0}(\Omega).
\end{eqnarray}
We use standard definitions of function spaces associated with $Q_{T}$; namely, $L^{q,r}(Q_T)$ consists of functions from $L^{r}(Q_T)$, $r \ge 1$, with the finite norm
\[
{\displaystyle \LN{g}_{q,r,Q_T}
:= \Big(\int_{0}^{T}\LN{g\left(\cdot,t\right)}_{q,\Omega}^{r}\: dt \Big)^{1/r}<\infty}.
\]
We denote by $H^{1}\left(Q_{T}\right)$ the space $L^{2}\left(0,T;H^{1}\left(\Omega\right)\right)$, and by $H_{0}^{1}\left(Q_{T}\right)$ the subspace of $H^{1}\left(Q_{T}\right)$ that consists of functions vanishing on $S_{T}$. The subspace $\mV \left(Q_{T}\right)$ consists of functions with the finite norm
\[\LN{w}_{\mV}^{2}:= \mbox{ess} \sup_{t \in \left(0,T\right)}
\LN{w \left(\cdot,t\right)}_{\Omega}^{2}+\LN{\nabla w}_{2,Q_{T}}^{2}.\]
Functions in $\mV \left(Q_{T}\right)$ have bounded $L^{2}$ norms at cross-sections of the time-domain. The space
\[\mV^{1,0} \left(Q_{T}\right):= H^{1} \left(Q_{T}\right) \cap C\left(0,T;L^{2} \left(\Omega\right)\right)\]
is a subspace of $\mV \left(Q_{T}\right)$. We denote by $\mV_{0}^{1,0}\left(Q_{T}\right)$ the subspace of $\mV^{1,0}$ that consists of functions vanishing on $S_T$ .

A function $u\in\mV_{0}^{1,0}\left(Q_{T}\right)$ is called the generalized solution of the problem (\ref{eq:Prob})-(\ref{eq:Cond_Bdry}) if it satisfies the following relation for all $w\in H_{0}^{1}\left(Q_{T}\right)$
\begin{align}
\int_{Q_{T}}\nabla u\cdot\nabla w \,dx \,dt
- \int_{Q_{T}}uw_{t} \,dx \,dt
+ \int_{Q_{T}} a \cdot \nabla u \,w \,dx \,dt & \label{eq:GenSol} \\
+ \int_{\Omega} \left(u\left(x,T\right) w \left(x,T\right)
- u \left(x,0\right) w \left(x,0\right)\right)\,dx
& = \int_{Q_{T}}fw \,dx \,dt. \nonumber
\end{align}
The existence of the generalized solution follows from well known results (see, e.g., \cite{Ladyzhen-85,LadyzhenSU-67}).

\section{Guaranteed Error Bounds for Space-conforming Approximations}

Let $v \in H_{0}^{1} \left(Q_{T}\right)$. In order to deduce a computable measure of $u-v$, we insert $v$ into the integral identity and represent it as follows:
\begin{alignat}{1}
\label{eq:IbPv_0}
& \int_{Q_{T}} \nabla \left(u - v\right) \cdot \nabla w \,dx \,dt
+ \int_{Q_{T}} \left(a \cdot \nabla \left(u - v\right)\right) w \,dx \,dt
- \underbrace{\int_{Q_{T}} \left(u - v\right) w_{t} \,dx \,dt}_{\left(*\right)} \\
& + \underbrace{\int_{\Omega} \left(\left(u \left(x,T\right)
- v\left(x,T\right)\right) w \left(x,T\right) - \left(u \left(x,0\right)
- v\left(x,0\right)\right) w \left(x,0\right)\right) \,dx}_{\left(**\right)} \nonumber \\
= & \int_{Q_{T}} \left(fw - \nabla v \cdot \nabla w - v_{t} w - \left(a \cdot \nabla v\right) w \right) \,dx \,dt. \nonumber
\end{alignat}
We note that
\begin{alignat}{1}
\frac{1}{2}\int_{\Omega} \left. \left| w \left(x,t\right)\right|^{2}\right|_{t=0}^{t=T} \,dx
& = \int_{Q_{T}} \frac{1}{2} \frac{d}{dt} \left(w^{2} \left(x,t\right)\right) \,dx \,dt
= \int_{Q_{T}}ww_{t} \,dx \,dt.\label{eq:IntDiff_Time}
\end{alignat}
Using (\ref{eq:IntDiff_Time}) and setting $w=u-v$, from the two integrals ({*}) and ({*}{*}) of (\ref{eq:IbPv_0}) we find
\begin{alignat}{1}
& - \frac{1}{2}\int_{\Omega} \left. \left|u \left(x,t\right) - v \left(x,t\right)\right|^{2}\right|_{t=0}^{t=T} \,dx
+ \int_{\Omega} \left. \left|u\left(x,t\right) - v \left(x,t\right)\right|^{2}\right|_{t=0}^{t=T} \,dx \label{eq:Int_Time}\\
= & \frac{1}{2}\int_{\Omega}\left.\left|u\left(x,t\right) - v\left(x,t\right)\right|^{2}\right|_{t=0}^{t=T} \,dx. \nonumber
\end{alignat}
Hence, we arrive at the relation
\begin{alignat}{1}
& \int_{Q_{T}} \left| \nabla \left(u - v\right)\right|^{2} \,dx \,dt
+ \int_{Q_{T}} \left(a \cdot \nabla \left(u - v\right)\right) \left(u - v\right) \,dx \,dt
+\left. \frac{1}{2}\LN{u - v}_{\Omega}^{2}\right|_{t=0}^{t=T} \label{eq:IbPv_1}\\
= & \int_{Q_{T}} \left(f \left(u - v\right) - \nabla v \cdot \nabla \left(u - v\right)
- v_{t}\left(u - v\right) - \left(a \cdot \nabla v\right) \left(u - v\right)\right) \,dx \,dt.\nonumber
\end{alignat}
Now, for $w = w\left(x,t\right)$ we use the identity
\begin{alignat}{1}
\Div\left(wa\right) & =w\:\Div a+a\cdot\nabla w,\label{eq:IdenDiv}
\end{alignat}
which holds for $t\in\left[0,T\right]$ almost everywhere. We integrate it over $Q_{T}$ and take into account that $w$ vanishes over its boundary $S_{T}$. We have
\[
\int_{Q_{T}} \left(a \cdot \nabla w \right) w \,dx \,dt = \int_{Q_{T}} wa \cdot \nabla w \,dx \,dt
= - \int_{Q_{T}} \Div \left(wa\right)w \,dx \,dt.
\]
Then, using (\ref{eq:IdenDiv}) in the last term, we find
\[
\int_{Q_{T}} \left(a \cdot \nabla w\right)w \,dx \,dt
= - \int_{Q_{T}} \left(w^{2} \Div a + wa \cdot \nabla w\right) \,dx \,dt.
\]
Hence
\begin{alignat}{1}
\int_{Q_{T}} \left(a \cdot \nabla w\right)w \,dx \,dt
& = - \frac{1}{2} \int_{Q_{T}}w^{2} \Div a \,dx \,dt .\label{eq:DivRel}
\end{alignat}
Using (\ref{eq:DivRel}), we rearrange (\ref{eq:IbPv_1}) as follows
\begin{alignat}{1}
& \int_{Q_{T}} \left| \nabla \left(u - v\right) \right|^{2} \,dx \,dt
- \frac{1}{2} \int_{Q_{T}}\Div a \, \left(u - v\right)^{2} \,dx \,dt
+ \left. \frac{1}{2} \LN{u - v}_{\Omega}^{2}\right|_{t=T} \label{eq:IbPv_2}\\
= & \int_{Q_{T}} \left(f - v_{t} - a \cdot \nabla v\right) \left(u - v\right) \,dx \,dt
- \int_{Q_{T}} \nabla v \cdot \nabla \left(u - v\right) \,dx \,dt
+ \left. \frac{1}{2}\LN{u-v}_{\Omega}^{2}\right|_{t=0}.\nonumber
\end{alignat}
Let us define the error $e=u-v$. Since $e=0$ on $\partial\Omega$, we have another relation for $t\in\left[0,T\right]$ almost everywhere
\begin{equation}
\int_{\Omega} \left(e \, \Div y + y \cdot \nabla e \right) \,dx
= \int_{\partial \Omega} e\, y \cdot \gamma \,dx = 0,
\label{eq:IdenDiv_y}
\end{equation}
where $\gamma$ denote the unit outward normal vector to $\partial\Omega$, and the vector-valued function $y=y\left(x,t\right)$ belongs to the space
\[
\Hdiv{Q_{T}}:= \left\{ y \in L^{2} \left(0,T;L^{2} \left(\Omega,\R^{d}\right)\right),\,
\Div y \in L^{2} \left(0,T;L^{2} \left(\Omega\right)\right)\right\}.
\]
Then, using \eqref{eq:IdenDiv_y} in (\ref{eq:IbPv_2}) we get
\begin{alignat}{1}
& \int_{Q_{T}} \Big(\left| \nabla e\right|^{2}
- \frac{1}{2} \Div a\, \left| e\right|^{2}\Big) \,dx \,dt
+ \frac{1}{2} \LN{e \left(\cdot,T\right)}_{\Omega}^{2}\label{eq:IbPv_3}\\
= & \int_{Q_{T}} r e \,dx \,dt
+\int_{Q_{T}}\left(y - \nabla v\right)\cdot\nabla e \,dx \,dt
+ \frac{1}{2}\LN{e \left(\cdot,0\right)}_{\Omega}^{2},\nonumber
\end{alignat}
where
\begin{equation}
\label{eq:Def_r}
r = r\left(v,y\right) := f - v_{t} - a \cdot \nabla v + \Div y.
\end{equation}
Henceforth, it is convenient to consider two cases; namely, $\Div a=0$, and $\Div a<0$.

\subsection{Case $\Div a = 0$}\label{sec:div0}

If $\Div a=0$ (which is typical if convection is defined by a constant vector $a$) then the left hand side (LHS) of (\ref{eq:IbPv_3}) contains the combined error norm
\begin{equation}
\CN{e}^{2} := \LN{\nabla e}_{Q_{T}}^{2}
+ \frac{1}{2} \LN{e \left(\cdot,T\right)}_{\Omega}^{2}.\label{eq:CN}
\end{equation}
We use the H\"older estimate
\begin{equation}
\int_{Q_{T}} r e \,dx \,dt
\le \LN{r}_{Q_{T}} \LN{e}_{Q_{T}}. \label{eq:Hoelder}
\end{equation}
Since the following relation holds for $t \in \left(0,T\right)$ almost everywhere
\[
\LN{e}_{\Omega} \le C_{F_{\Omega}} \LN{\nabla e}_{\Omega},
\]
we have
\begin{equation}
\LN{e}_{Q_{T}}
\le C_{F_{\Omega}} \LN{\nabla e}_{Q_{T}},\label{eq:Friedrich}
\end{equation}
where $C_{F_{\Omega}}$ is the constant in the Friedrichs inequality. Therefore,
\begin{equation}
\int_{Q_{T}} r e \,dx \,dt
\le C_{F_{\Omega}} \LN{r}_{Q_{T}}\LN{\nabla e}_{Q_{T}}
\le C_{F_{\Omega}}\LN{r}_{Q_{T}}\CN{e}.\label{eq:Friedrich_y}
\end{equation}
Moreover,
\begin{equation}
\int_{Q_{T}} \left(y - \nabla v \right) \cdot \nabla e \,dx \,dt
\le\LN{y - \nabla v}_{Q_{T}} \LN{\nabla e}_{Q_{T}} \le \LN{y - \nabla v}_{Q_{T}}\CN{e}.\label{eq:Schwarz_y}
\end{equation}
Now (\ref{eq:IbPv_3}), (\ref{eq:Friedrich_y}), and (\ref{eq:Schwarz_y}) imply
\[
\CN{e}^{2}
\le \LN{y-\nabla v}_{Q_{T}}\CN{e} + C_{F_{\Omega}} \LN{r}_{Q_{T}} \CN{e}
+ \frac{1}{2}\LN{e \left(\cdot,0\right)}_{\Omega}^{2}.
\]
Setting
\begin{equation}
\mM_{Q_{T}}\left(v, y\right) := \LN{y - \nabla v}_{Q_{T}} + C_{F_{\Omega}} \LN{r}_{Q_{T}},
\label{eq:MajSpace}
\end{equation}
and using some simple arithmetic-geometric inequalities, we deduce the simplest form of the \emph{guaranteed} upper bound
\begin{alignat}{1}
\label{eq:UB_div0}
2 \CN{e} & \le \mM_{Q_{T}}\left(v, y\right) + \Big(\mM_{Q_{T}}^{2}\left(v, y\right) +2 \LN{e \left(\cdot,0\right)}_{\Omega}^{2}\Big)^{1/2}.
\end{alignat}
It is easy to see that if the initial conditions are exactly satisfied, i.e.
\[v\left(x,0\right)=\varphi\left(x\right),\]
then we get a simple form of the estimate
\begin{alignat}{1}
\CN{e} & \le \mM_{Q_{T}}\left(v, y\right), \label{eq:UB_div0_1}
\end{alignat}
which reflects the fact that the combined error norm is controlled by integral norms of discrepancies in the basic relations
\begin{alignat*}{2}
y - \nabla v & = 0 & \quad & \mbox{in }Q_{T}\,,\\
\Div y + f - v_{t} - a \cdot \nabla v & = 0 & \quad & \mbox{in }Q_{T}\,.
\end{alignat*}
It is also possible to bound the term $\LN{e\left(\cdot, T\right)}_{\Omega}$ alone. Using Young inequality in (\ref{eq:Friedrich_y}) we get
\begin{equation}
\int_{Q_{T}} r e \,dx \,dt
\le \int_{0}^{T} \frac{C_{F_{\Omega}}^{2}}{2 \alpha} \LN{r}_{\Omega}^{2} \,dt
+ \int_{0}^{T} \frac{\alpha}{2} \LN{\nabla e}_{\Omega}^{2} \,dt,
\label{eq:FriedrichYoung_y}
\end{equation}
where $\alpha \left(t\right)$ is some positive bounded function, i.e.
\[0 < \alpha \left(t\right) \in L^{\infty} \left(0,T\right).\]
Using Young inequality in (\ref{eq:Schwarz_y}) with some
\[0 < \beta \left(t\right) \in L^{\infty} \left(0,T\right),\]
we get
\begin{equation}
\int_{Q_{T}} \left(y - \nabla v \right) \cdot \nabla e \,dx \,dt
\le \int_{0}^{T} \frac{1}{2 \beta} \LN{y - \nabla v}_{\Omega}^{2} \,dt
+ \int_{0}^{T} \frac{\beta}{2} \LN{\nabla e}_{\Omega}^{2} \,dt.
\label{eq:SchwarzYoung_y}
\end{equation}
Then, (\ref{eq:IbPv_3}), (\ref{eq:FriedrichYoung_y}), and (\ref{eq:SchwarzYoung_y}) imply
\begin{alignat}{1}
& \int_{0}^{T} \Big(1 - \frac{\alpha}{2} - \frac{\beta}{2}\Big) \LN{\nabla e}_{\Omega}^{2} \,dt
+ \frac{1}{2}\LN{e\left(\cdot,T\right)}_{\Omega}^{2} \\
\leq & \int_{0}^{T} \frac{C_{F_{\Omega}}^{2}}{2 \alpha} \LN{r}_{\Omega}^{2} \,dt
+ \int_{0}^{T} \frac{1}{2 \beta} \LN{y - \nabla v}_{\Omega}^{2} \,dt
+ \frac{1}{2}\LN{e\left(\cdot, 0\right)}_{\Omega}^{2}. \nonumber
\end{alignat}
Choosing $\alpha = \beta = 1$ we finally get
\begin{align}
\label{eq:eTest_div0}
\LN{e\left(\cdot,T\right)}_{\Omega}^{2}
& \leq \mM_{T}^{2}\left(v,y\right),
\end{align}
where
\begin{align*}
\mM_{T}^{2}\left(v,y\right)
& := C_{F_{\Omega}}^{2} \LN{r}_{Q_{T}}^{2}
+ \LN{y - \nabla v}_{Q_{T}}^{2}
+ \LN{e\left(\cdot, 0\right)}_{\Omega}^{2}. \nonumber
\end{align*}

\subsection{Case $\Div a < 0$}\label{sec:delta}

If ${\displaystyle - \frac{1}{2}\Div a = \delta^{2}>0}$, then we obtain estimates in terms of a different weighted norm
\begin{equation}
\label{eq:CN_delta}
\CN{e}_{\delta}^{2}
:= \LN{\nabla e}_{Q_{T}}^{2} + \LN{\delta e}_{Q_{T}}^{2}
+ \LN{e\left(\cdot,T\right)}_{\Omega}^{2}.
\end{equation}
By (\ref{eq:IbPv_3}) we find that
\begin{alignat}{1}
& \LN{\nabla e}_{Q_{T}}^{2} + \LN{\delta e}_{Q_{T}}^{2}
+ \frac{1}{2} \LN{e \left(\cdot,T\right)}_{\Omega}^{2} \label{eq:UB_delta}\\
\le & \int_{Q_{T}} r e \,dx \,dt
+ \int_{0}^{T} \LN{y - \nabla v}_{\Omega}\LN{\nabla e}_{\Omega} \,dt
+ \frac{1}{2} \LN{e \left(\cdot,0\right)}_{\Omega}^{2}. \nonumber
\end{alignat}
For reasons that will become clear later, we introduce a function $\lambda\left(x,t\right)$ with values in $\left[0,1\right]$, and rewrite the first term on the right hand side (RHS) of (\ref{eq:UB_delta}) as follows
\[
\int_{Q_{T}} r e \,dx \,dt
= \int_{Q_{T}}\lambda \, re \,dx \,dt + \int_{Q_{T}} \left(1 - \lambda\right)re \,dx \,dt.
\]
Then, for some $0 < \alpha \left(t\right) \in L^{\infty} \left(0,T\right)$ we have
\begin{alignat}{1}
\label{eq:Est_M2}
\int_{Q_{T}} r e \,dx \,dt
\le & \int_{0}^{T} \Big(\frac{1}{4 \alpha} \LN{\frac{\lambda}{\delta}r}_{\Omega}^{2}
+ \alpha\LN{\delta e}_{\Omega}^{2}\Big) \,dt
+ \int_{0}^{T} \LN{\left(1 - \lambda \right) r}_{\Omega} \LN{e}_{\Omega} \,dt \\
\le & \int_{0}^{T} \Big(\frac{1}{4 \alpha} \LN{\frac{\lambda}{\delta}r}_{\Omega}^{2}
+ \alpha \LN{\delta e}_{\Omega}^{2}\Big) \,dt
+ C_{F_{\Omega}} \int_{0}^{T}\LN{\left(1 - \lambda \right)r}_{\Omega}
\LN{\nabla e}_{\Omega} \,dt.\nonumber
\end{alignat}
Using the Young inequality in the last term of (\ref{eq:Est_M2}) with some
\[0 < \gamma \left(t\right) \in L^{\infty} \left(0,T\right),\]
we get
\begin{align}
\label{eq:Est_M2_1}
\int_{Q_{T}} r e \,dx \,dt
\le & \int_{0}^{T} \Big(\frac{1}{4 \alpha} \LN{\frac{\lambda}{\delta}r}_{\Omega}^{2}
+ \alpha \LN{\delta e}_{\Omega}^{2}\Big) \,dt
+ \int_{0}^{T} \Big( \frac{C_{F_{\Omega}}^{2}}{2 \gamma} \LN{\left(1 - \lambda \right) r}_{\Omega}^{2}
+ \frac{\gamma}{2} \LN{\nabla e}_{\Omega}^{2} \Big) \,dt.
\end{align}
Similarly, for some $0 < \beta \left(t\right) \in L^{\infty} \left(0,T\right)$, we have
\begin{alignat}{1}
\int_{0}^{T}\LN{y-\nabla v}_{\Omega}\LN{\nabla e}_{\Omega} \,dt
& \le \int_{0}^{T}\frac{\beta}{2}\LN{\nabla e}_{\Omega}^{2} \,dt
+ \int_{0}^{T}\frac{1}{2\beta} \LN{y - \nabla v}_{\Omega}^{2} \,dt.\label{eq:Est_M1}
\end{alignat}
By (\ref{eq:UB_delta}-\ref{eq:Est_M1}), we conclude that
\begin{alignat}{1}
\label{mainest}
& \int_{0}^{T} \big(1 - \frac{\beta}{2} - \frac{\gamma}{2}\big) \LN{\nabla e}_{\Omega}^{2} \,dt
+ \int_{0}^{T}\left(1 - \alpha\right) \LN{\delta e}_{\Omega}^{2} \,dt
+ \frac{1}{2} \LN{e \left(\cdot,T\right)}_{\Omega}^{2} \\
\le & \int_{0}^{T} \frac{1}{4 \alpha} \LN{\frac{\lambda}{\delta} r}_{\Omega}^{2} \,dt
+ \int_{0}^{T} \frac{C_{F_{\Omega}}^{2}}{2\gamma} \LN{\left(1 - \lambda\right)r}_{\Omega}^{2} \,dt
+ \int_{0}^{T}\frac{1}{2\beta}\LN{y - \nabla v}_{\Omega}^{2} \,dt
+ \frac{1}{2} \LN{e \left(\cdot,0\right)}_{\Omega}^{2}, \nonumber
\end{alignat}
where $\alpha(t)$, $\beta(t)$, and $\gamma(t)$ are some arbitrary positive functions satisfying the conditions
\begin{equation}
2 - \beta - \gamma > 0, \qquad 1 - \alpha > 0.\label{eq:CondFunc}
\end{equation}
%

In particular, if we set $\alpha=\beta=\gamma=1$ in (\ref{mainest}), we can bound the term $\LN{e\left(\cdot,T\right)}_{\Omega}$ alone as
\begin{alignat}{1}
\label{eq:eTest_delta}
\LN{e\left(\cdot,T\right)}_{\Omega}^{2}
\le & \, \mN_{T}^{2}\left(\lambda,v,y\right),
\end{alignat}
where
\begin{align*}
\mN_{T}^{2}\left(\lambda,v,y\right)
:= & \, \frac{1}{2} \LN{\frac{\lambda}{\delta}r}_{Q_T}^{2}
+ C_{F_{\Omega}}^{2}\LN{\left(1 - \lambda\right)r}_{Q_T}^{2}
+ \LN{y - \nabla v}_{Q_T}^{2}
+ \LN{e \left(\cdot,0\right)}_{\Omega}^{2}.\nonumber
\end{align*}
%
If we set ${\displaystyle \alpha = \beta = \gamma = \frac{1}{2}}$, then we obtain the estimate in terms of the weighted space-time norm
\begin{alignat}{1}
\label{eq:UB_delta_2}
\CN{e}_{\delta}^{2}
& \le 2 \mN_{T}^{2}\left(\lambda,v,y\right) - \LN{e \left(\cdot,0\right)}_{\Omega}^{2}.
\end{alignat}
\begin{rem}
The function $\lambda$ involved in these estimates can be useful if $\delta$ attains small values at some points of $Q_{T}$. In this case, $\frac{1}{\delta}$ is a large penalty of the term $r\left(v,y\right)$, and the estimate may be too pessimistic unless the value of $r\left(v,y\right)$ is very small. In such a case, the function $\lambda\left(x,t\right)$ can be used to compensate this drawback. Indeed, we a priori know $a$ and $\Div a$. Therefore, we can select $\lambda$ in such a way that $\frac{\lambda}{\delta}$ remains of order $1$. If $\delta$ is sufficiently large everywhere, then it may be useful to simply set $\lambda=1$ and use the simplified estimate
\begin{alignat}{1}
& \int_{0}^{T} \Big(1 - \frac{\beta}{2}\Big) \LN{\nabla e}_{\Omega}^{2} \,dt
+ \int_{0}^{T}\left(1 - \alpha \right)\LN{\delta e}_{\Omega}^{2} \,dt
+ \frac{1}{2} \LN{e \left(\cdot,T\right)}_{\Omega}^{2} \label{eq:UB_delta_3}\\
\le & \int_{0}^{T}\frac{1}{4\alpha}\LN{\frac{1}{\delta}r}_{\Omega}^{2} \,dt
+ \int_{0}^{T}\frac{1}{2\beta}\LN{y-\nabla v}_{\Omega}^{2} \,dt
+ \frac{1}{2} \LN{e \left(\cdot,0\right)}_{\Omega}^{2} .\nonumber
\end{alignat}
Further, when the initial conditions are exactly satisfied, i.e.
\[v \left(x,0\right) = \varphi \left(x\right),\]
by choosing $\alpha=1/2,$ and $\beta=1$, we have a simpler estimate
\begin{alignat}{1}
\CN{e}_{\delta}^{2}
& \le \LN{\frac{1}{\delta}r}_{Q_{T}}^{2}
+ \LN{y - \nabla v}_{Q_{T}}^{2}.\label{eq:UB_delta_4}
\end{alignat}
\end{rem}

\section{Nonconforming approximations}

\subsection{Method 1}
Assume that we have obtained some nonconforming approximation $\hat{v}$. The simplest approach to control its accuracy is to project it on a certain space of conforming approximations (e.g., by using averaging or other post-processing techniques). For elliptic type problems this method has been thoroughly discussed and tested in \cite{LazarovRT-09}. Let $P$ be some suitable mapping such that
\[P \hat{v} \in H^{1}_{0}(Q_{T}).\]
Then we can apply (\ref{mainest}) (or (\ref{eq:UB_delta_2}) to this function. By the triangle inequality we obtain
\begin{equation}
\label{nonconf1}
\CN{u - \hat{v}}_{\delta}
\leq \CN{\hat{v} - P \hat{v}}_{\delta} + \CN{u - P \hat{v}}_{\delta},
\end{equation}
where the first term on the RHS is directly computable and represents the \textit{nonconformity error}, whereas the second term can be estimated by the estimates derived in the previous section. It is worth noting that, in this approach, we do not exploit any specific structure of the underlying discretization method, and therefore, this approach is valid for any nonconforming approximation.

\subsection{Method 2}
We consider a special (but practically valuable) case of incremental approximations, in which the corresponding function $\hat v$ is uniformly bounded and piecewise  smooth (more precisely, it is smooth with respect to spatial and time variables on each time interval $(t_k,t_{k+1})$) and may have jumps in time at points $t_k$, $k=0,1,...N$, where $t_0=0$ and $t_N=T$. We show that the corresponding error estimate can be directly derived by a certain limit procedure applied to the main integral relation ({\ref{eq:IbPv_2}}).
Indeed, from (\ref{eq:IbPv_2}) we see that only the term
\[-\int_{Q_{T}}v_{t}\left(u-v\right) dx \, dt\]
involves the time derivative. To outline how the procedure acts in the context of time-nonconformity, we first consider the case where $\hat v$ has only one jump at a point $t = \tau \in (0,T)$.
Assume that
\[\hat v(\tau-0)=v^-, \quad \hat v(\tau+0)=v^+, \quad v^-\not= v^+ \,.\]
On $(0,T)$ we construct a sequence of functions  $v^\epsilon$ ($\epsilon$ is a small positive number) as follows:
\[
\begin{cases}
v^{\epsilon}(x,t)=\hat v(x,t), & \quad t \in (0,\tau - \epsilon) \cup (\tau, T),\\
v^\epsilon(x,t)=\hat v(x,\tau-\epsilon)+(v^+-v^-)\,\frac{\left( t - \tau + \epsilon \right)}{\epsilon},
& \quad t \in (\tau - \epsilon, \tau)\\
\end{cases}.
\]
It is easy to see that $v^\epsilon(x,t)$ tends to $\hat{v}\left(x,t\right)$ in $L^{2}(Q_T)$ as $\epsilon\rightarrow 0$. Moreover,
\[
\nabla v^\epsilon(x,t)
=\left\{\begin{array}{ll}
\nabla\hat v(x,\tau),\; & t \in (0,\tau-\epsilon),\\
\nabla\hat v(x,\tau-\epsilon)+(\nabla v^+-\nabla v^-)\,\frac{\left( t - \tau + \epsilon \right)}{\epsilon},\;
& t \in (\tau-\epsilon,\tau)
\end{array}\right.\,,
\]
and the space gradients also converge to $\nabla\hat v(x,\tau)$ in $L^{2}$.
We now apply the relation (\ref{eq:IbPv_2}) to $v^\epsilon$:
\begin{alignat}{1}
& \int_{Q_{T}}\left|\nabla\left(u-v^\epsilon\right)\right|^{2} \,dx \,dt
- \frac{1}{2}\int_{Q_{T}} \Div a\:\left(u-v^\epsilon\right)^{2} \,dx \,dt
+ \left.\frac{1}{2} \LN{u-\hat v}_{\Omega}^{2}\right|_{t=T}\label{eq:veps}\\
= & \int_{Q_{T}}\left(f-v^\epsilon_{t}-a\cdot\nabla v^\epsilon\right)\left(u-v^\epsilon\right) \,dx \,dt
- \int_{Q_{T}} \nabla v^\epsilon\cdot\nabla \left(u-v^\epsilon\right) \,dx \,dt
+ \left.\frac{1}{2}\LN{u-\hat v}_{\Omega}^{2}\right|_{t=0}.\nonumber
\end{alignat}
Note that, as $\epsilon\rightarrow\,0$, the LHS of (\ref{eq:veps}) tends to
\[
\int_{Q_{T}}\left|\nabla\left(u-\hat v\right)\right|^{2} \,dx \,dt
- \frac{1}{2}\int_{Q_{T}} \Div a\:\left(u-\hat v\right)^{2} \,dx \,dt
+ \left.\frac{1}{2} \LN{u-\hat v}_{\Omega}^{2}\right|_{t=T}.
\]
Moreover,
\begin{align*}
\int_{Q_{T}}
\nabla v^\epsilon\cdot\nabla \left(u-v^\epsilon\right) \,dx \,dt
& \rightarrow\,\int_{Q_{T}} \nabla \hat v\cdot\nabla \left(u-\hat v\right) \,dx \,dt,
\end{align*}
and
\begin{align*}
\int_{Q_{T}}\left(f-a\cdot\nabla v^\epsilon\right)\left(u-v^\epsilon\right) \,dx \,dt\,
& \rightarrow\, \int_{Q_{T}}\left(f-a\cdot\nabla \hat v\right)\left(u-\hat v\right) \,dx \,dt.
\end{align*}
It remains to consider the term involving the time derivative. We have
\begin{alignat}{1}
\IntO\int_{0}^{\tau}v_{t}^{\epsilon}u \,dx \,dt
& \rightarrow\IntO\int_{0}^{\tau-\epsilon}\hat{v}_{t}u \,dx \,dt
+ \IntO\int_{\tau-\epsilon}^{\tau}v_{t}^{\epsilon}u \,dx \,dt
\label{eq:NC_vtu} \\
& \rightarrow \IntO \int_{0}^{\tau-\epsilon}\hat{v}_{t}u \,dx \,dt
+ \IntO \int_{\tau-\epsilon}^{\tau} \left( \hat{v}_{t}
+ \frac{v^{+} - v^{-}}{\epsilon} \right) u \,dx \,dt \nonumber \\
& \rightarrow \IntO \int_{0}^{\tau}\hat{v}_{t}u \,dx \,dt
+ \IntO \frac{v^{+} - v^{-}}{\epsilon} \int_{\tau-\epsilon}^{\tau} u \,dx \,dt \nonumber \\
& \rightarrow\IntO\int_{0}^{\tau}\hat{v}_{t}u \,dx \,dt
+\IntO \left(v^{+} - v^{-}\right)u\left(\tau\right)\,dx, \nonumber
\end{alignat}
since $u\left(t\right)$ is continuous at $t=\tau$. Moreover,
\begin{alignat}{1}
\IntO\int_{0}^{\tau}v_{t}^{\epsilon}v^{\epsilon} \,dx \,dt
& = \IntO\int_{0}^{\tau-\epsilon}v_{t}^{\epsilon}v^{\epsilon} \,dx \,dt
+ \IntO\int_{\tau-\epsilon}^{\tau}v_{t}^{\epsilon}v^{\epsilon} \,dx \,dt \label{eq:NC_vtv} \\
& = \IntO\int_{0}^{\tau-\epsilon}\hat v_{t}\hat v \,dx \,dt
+ \IntO\int_{\tau-\epsilon}^{\tau}  \frac{1}{2}\frac{d}{dt} \left(v^{\epsilon}\right)^{2} \,dx \,dt \nonumber \\
& \rightarrow \IntO\int_{0}^{\tau}\hat v_{t}\hat v \,dx \,dt
+\frac{1}{2}\IntO\left(\left(v^{+}\right)^{2} - \left(v^-\right)^{2}\right) \,dx. \nonumber
\end{alignat}
From (\ref{eq:NC_vtu}) and (\ref{eq:NC_vtv}) it follows that
\begin{alignat*}{1}
& - \int_{Q_{T}}v^\epsilon_{t}\left(u-v^\epsilon\right) \,dx \,dt \\
\rightarrow & -\int_{Q_{T}}\hat{v}_{t} \left(u-\hat{v}\right) \,dx \,dt
- \int_{\Omega}\big(v^{+}-v^-\big) \Big(u(\tau) - \frac{v^{+}+v^-}{2}\Big) \,dx \\
= & - \int_{Q_{T}}\hat{v}_{t} \left(u-\hat{v}\right)
- \int_{\Omega}\left(v^{+}-v^-\right) \left(u\left(\tau\right)-v^-\right) \,dx
+ \frac{1}{2}\int_{\Omega}\left(v^{+}-v^-\right)^{2} \,dx,
\end{alignat*}
and we find that
\begin{alignat}{1}
\label{eq:mainest_NC}
& \int_{Q_{T}}\left|\nabla \hat{e} \right|^{2} \,dx \,dt
- \frac{1}{2}\int_{Q_{T}} \Div a \,(\hat{e})^{2} \,dx \,dt
+\frac{1}{2} \LN{\hat{e}\left(\cdot, T \right)}_{\Omega}^{2} \\
= & \int_{Q_{T}}\left(f - \hat{v}_t - a \cdot \nabla \hat v\right) \hat{e} \,dx \,dt
+\int_{Q_{T}} \nabla \hat v \cdot \nabla \hat{e} \,dx \,dt \nonumber\\
& - \int_{\Omega}\left(v^{+}-v^-\right) \left(u\left(\tau\right)-v^-\right) \,dx
+ \frac{1}{2}\int_{\Omega}\left(v^{+} - v^{-}\right)^{2} \,dx
+ \frac{1}{2}\LN{\hat{e}\left(\cdot,0\right)}_{\Omega}^{2} \nonumber \\
\leq \, & C_{F_{\Omega}}\LN{\hat{r}}_{Q_{T}} \LN{\nabla \hat{e}}_{Q_{T}}
+ \LN{y-\nabla\hat{v}}_{Q_{T}} \LN{\nabla \hat{e}}_{Q_{T}} \nonumber \\
& + \LN{v^{+} - v^{-}}_{\Omega} \LN{u(\tau) - v^{-}}_{\Omega}
+ \frac{1}{2} \LN{v^{+} - v^{-}}^2_{\Omega}
+ \frac{1}{2} \LN{\hat{e}\left(\cdot,0\right)}_{\Omega}^{2}, \nonumber
\end{alignat}
where $\hat{r} = r\left(\hat{v},y\right)$, and $\hat{e} = u - \hat{v}$. Now recall that for $u(\tau) - v^{-}$ we can apply the majorant (\ref{eq:eTest_div0}) or (\ref{eq:eTest_delta}) depending on whether $\Div a=0$ or $\Div a <0$, respectively. We proceed with both cases separately.

\subsubsection{\textbf{Case $\Div a = 0$}}
Note that $u(\tau) - v^{-} = e\left(\cdot, \tau \right)$. Therefore, using (\ref{eq:eTest_div0}) for $u(\tau) - v^{-}$ we get
\begin{equation}
\label{eq:etau_est_div0}
\LN{u(\tau) - v^{-}}_{\Omega}
= \LN{e\left(\cdot, \tau \right)}_{\Omega}
\leq \mM_{\tau}\left(\hat{v},y\right).
\end{equation}
To simplify the notations, we introduce
\begin{equation}
\label{eq:CJ}
\mathcal{C}_{J}
= \LN{v^{+} - v^{-}}_{\Omega} \mM_{\tau}\left(\hat{v},y\right)
+ \frac{1}{2} \LN{v^{+} - v^{-}}^2_{\Omega}
+ \frac{1}{2} \LN{\hat{e}\left(\cdot,0\right)}_{\Omega}^{2}.
\end{equation}
Then, using (\ref{eq:etau_est_div0}) and (\ref{eq:CJ}) in (\ref{eq:mainest_NC}), and that $\Div a = 0$, we get
\begin{alignat}{1}
\label{eq:est_NC_1}
\LN{\nabla \hat{e}}_{Q_{T}}^{2}
+ \frac{1}{2} \LN{\hat{e}\left(\cdot, T \right)}_{\Omega}^{2}
\leq & \big(C_{F_{\Omega}}\LN{\hat{r}}_{Q_{T}}
+ \LN{y - \nabla\hat{v}}_{Q_{T}}\big) \LN{\nabla \hat{e}}_{Q_{T}}
+ \mathcal{C}_{J}.
\end{alignat}
Since $\LN{\nabla \hat{e}}_{Q_{T}} \le \CN{\hat{e}}$, we have
\begin{alignat}{1}
\CN{\hat{e}}^{2}
& \le \big(C_{F_{\Omega}}\LN{\hat{r}}_{Q_{T}}
+ \LN{y-\nabla\hat{v}}_{Q_{T}}\big) \CN{\hat{e}}
+ \mathcal{C}_{J}.
\end{alignat}
Using the definition \eqref{eq:MajSpace} and some simple arithmetic-geometric inequalities, we finally obtain
\begin{equation}
\label{eq:mainest_NC_div0}
2 \CN{\hat{e}} \le \mM_{Q_{T}}\left(\hat{v}, y\right)
+ \Big(\mM_{Q_{T}}^{2}\left(\hat{v}, y\right) + 4 \mathcal{C}_{J} \Big)^{1/2}.
\end{equation}
\subsubsection{\textbf{Case $\Div a < 0$}}
In this case we use (\ref{eq:eTest_delta}) for $u(\tau) - v^{-}$, which gives
\begin{equation}
\label{eq:etau_est_delta}
\LN{u(\tau) - v^{-}}_{\Omega}
= \LN{e\left(\cdot, \tau \right)}_{\Omega}
\leq \mN_{\tau}\left(\lambda_{1}, \hat{v}, y\right),
\end{equation}
where $\lambda_{1}\left(x,t\right)$ is a function with values in $[0,1]$ for $t \in \left(0,\tau\right)$. To simplify the notations, we now introduce
\begin{equation}
\label{eq:DJ}
\mathcal{D}_{J}
= \LN{v^{+} - v^{-}}_{\Omega} \mN_{\tau}\left(\lambda_{1}, \hat{v}, y\right)
+ \frac{1}{2} \LN{v^{+} - v^{-}}^2_{\Omega}
+ \frac{1}{2} \LN{\hat{e}\left(\cdot,0\right)}_{\Omega}^{2}.
\end{equation}
Proceeding in the same way as in the conforming case, we finally obtain
\begin{alignat}{1}
\label{mainest_NC_delta}
& \int_{0}^{T} \big(1 - \frac{\beta}{2} - \frac{\gamma}{2}\big) \LN{\nabla \hat{e}}_{\Omega}^{2} \,dt
+ \int_{0}^{T}\left(1 - \alpha\right) \LN{\delta \hat{e}}_{\Omega}^{2} \,dt
+ \frac{1}{2} \LN{\hat{e} \left(\cdot,T\right)}_{\Omega}^{2} \\
\le \,\, & \mathcal{D}_{J}
+ \int_{0}^{T} \frac{1}{4 \alpha} \LN{\frac{\lambda}{\delta} \hat{r}}_{\Omega}^{2} \,dt
+ \int_{0}^{T} \frac{C_{F_{\Omega}}^{2}}{2\gamma} \LN{\left(1 - \lambda\right) \hat{r}}_{\Omega}^{2} \,dt
+ \int_{0}^{T}\frac{1}{2 \beta} \LN{y - \nabla \hat{v}}_{\Omega}^{2} \,dt. \nonumber
\end{alignat}
If we set ${\displaystyle \alpha = \beta = \gamma = \frac{1}{2}}$, then we get a simpler estimate
\begin{equation}
\CN{\hat{e}}^{2}_{\delta}
\le 2 \mathcal{D}_{J}
+ 2 \mN_{T}^{2}\left(\lambda, \hat{v}, y\right) - 2 \LN{e \left(\cdot,0\right)}_{\Omega}^{2}.
\end{equation}
\begin{rem}
If, instead of (\ref{eq:CN_delta}), we introduce a slightly different weighted norm
\begin{equation}
\label{eq:CN_delta2}
\CN{e}_{\hat{\delta}}^{2}
:= \LN{\nabla e}_{Q_{T}}^{2} + \LN{\delta e}_{Q_{T}}^{2}
+ \frac{1}{2} \LN{e\left(\cdot,T\right)}_{\Omega}^{2},
\end{equation}
then, using (\ref{eq:etau_est_delta}) and (\ref{eq:DJ}) in (\ref{eq:mainest_NC}), and that $\LN{\nabla \hat{e}}_{Q_{T}} \le \CN{\hat{e}}_{\hat{\delta}}$, we get
\begin{align}
\CN{\hat{e}}_{\hat{\delta}}^{2}
& \le \big(C_{F_{\Omega}}\LN{\hat{r}}_{Q_{T}}
+ \LN{y-\nabla\hat{v}}_{Q_{T}}\big) \CN{\hat{e}}_{\hat{\delta}}
+ \mathcal{D}_{J},
\end{align}
which easily leads to the following estimate similar to (\ref{eq:mainest_NC_div0})
\begin{equation}
\label{eq:mainest_NC_delta}
2 \CN{\hat{e}}_{\hat{\delta}}
\le \mM_{Q_{T}}\left(\hat{v}, y\right)
+ \Big(\mM_{Q_{T}}^{2}\left(\hat{v}, y\right) + 4 \mathcal{D}_{J} \Big)^{1/2}.
\end{equation}
\end{rem}
We see that the estimates in both the cases (through (\ref{eq:CJ}) and (\ref{eq:DJ})) involve penalty terms depending on the jump $\LN{v^{+} - v^{-}}_{\Omega}$. If the latter quantity is small, then the overall value of the majorant does not essentially increase. If we wish to introduce more time-discontinuity points then we can easily extend these estimates using the techniques discussed to compute $\LN{e\left(\cdot, T \right)}_{\Omega}$ in Subsections~\ref{sec:div0} and \ref{sec:delta}.
\begin{rem}
It may be very convenient to use approximations discontinuous in time if the spatial discretizations are reconstructed during the process of time integration. In this case, at certain time moment $\tau$, we need to change the structure of the finite-dimensional space which is used for approximation of the spatial components of the solution. Then, we may have difficulties in conforming continuation of the approximate solution. With the help of nonconforming extensions (which are technically simple) we can easily obtain $v\left(\tau + 0 \right)$ by interpolating $v\left(\tau - 0 \right)$. Jumps of discontinuities that arises in such a procedure can be taken into account due to the penalty terms in the above estimates.
\end{rem}

\section*{Acknowledgment}
The authors are grateful to the Russian Foundation of Fundamental Sciences (grant No. 08-01-00655-a) and the Austrian Academy of Sciences for support.



\begin{thebibliography}{363}

\bibitem{GaevsRepin-05}
A.~Gaevskaya and S.~Repin: A posteriori error estimates for approximate solutions of linear parabolic problems,
\textit{Differ. Uravn.} \textbf{41 (7),} 970--983, 2005 (in Russian); translation in \textit{Differ. Equ.} \textbf{41 (7),} 970--983, 2005.
%
\bibitem{Ladyzhen-85}
O.A.~Ladyzhenskaya: \textit{The boundary value problems of mathematical physics,} Springer, New York, 1985.
%
\bibitem{LadyzhenSU-67}
O.A.~Ladyzhenskaya, V.A.~Solonnikov and N.N.~Uraltseva: \textit{Linear and Quasilinear Equations of Parabolic Type,} Nauka, Moscow, 1967.
%
\bibitem{LazarovRT-09}
R.~Lazarov, S.~Repin, and S.~Tomar: Functional a posteriori error estimates for discontinuous Galerkin method,
\textit{Numer. Methods Partial Differential Equations} \textbf{25,} 952--971, 2009.
%
\bibitem{NeittanRepin-10}
P.~Neittanam\"aki and S.~Repin: A posteriori error majorants for approximations of the evolutionary Stokes problem.
\textit{J. Numer. Math.} \textbf{18 (2),} 119--134, 2010.
%
\bibitem{Repin-02}
S.~Repin: Estimates of deviation from exact solutions of initial-boundary value problems for the heat equation,
\textit{Rend. Mat. Acc. Lincei} \textbf{13,} 121-133, 2002.
%
\bibitem{Repin-08}
S.~Repin: \textit{A posteriori estimates for partial differential equations,} Walter de Gruyter, Berlin, 2008.
%
\end{thebibliography}
\end{document}